\documentclass[11pt]{article}
\usepackage{latexsym}
\usepackage{amsmath, amsthm}
\usepackage[all,dvips,ps]{xy}
\usepackage{amsfonts}
\usepackage{mathrsfs}

\newtheorem{lem}{Lemma}[section]
\newtheorem{thm}{Theorem}[section]
\newtheorem{cor}{Corollary}[section]
\newtheorem{ass}{Assumption}[section]
\newtheorem{exa}{Example}[section]

\newtheorem{defn}{Definition}[section]

\newcommand{\Rs}{\mathbb{R}}

\newcommand{\Es}{\mathscr{E}}
\newcommand{\Sn}{{\cal S}_n }

\newcommand{\bz}{{\bf 0} }
\newcommand{\tl}{\langle}
\newcommand{\tr}{\rangle }

\newcommand{\rr}{\bar{r} }

\newcommand{\trace}{{\rm trace\,}}

\newcommand{\bpr}{{\bf Proof.} \hspace{1 em}}
\newcommand{\epr}{ \\ \hspace*{4.5in} $\Box$ }
\newcommand{\beq}{ \begin{equation} }
\newcommand{\eeq}{ \end{equation} }
\newcommand{\bt}{ \begin{tabular} }
\newcommand{\et}{ \end{tabular} }

\def\rb{\bar{r}}
\def\Zh{\hat{Z}}
\def\zh{\hat{z}}

\begin{document}

\bibliographystyle{plain}
\title{On Affine Motions and Universal Rigidity of \\ Tensegrity Frameworks
 \thanks{Research supported by the Natural Sciences and Engineering
         Research Council of Canada.} }
\vspace{0.3in}
        \author{ A. Y. Alfakih
  \thanks{E-mail: alfakih@uwindsor.ca}
  \\
          Department of Mathematics and Statistics \\
          University of Windsor \\
          Windsor, Ontario N9B 3P4 \\
          Canada
\and
               Viet-Hang Nguyen
\thanks{Work done while visiting the University of Windsor.
E-mail:Viet-Hang.Nguyen@grenoble-inp.fr}
\\
         Laboratoire G-SCOP, Grenoble UJF \\
         46 avenue F\'{e}lix Viallet  \\
         38000 Grenoble \\
         France }

\date{\today}
\maketitle


\vspace{0.1in}

\begin{abstract}
Recently, Alfakih and Ye [Lin. Algebra Appl. 438:31--36, 2013] 
proved that if an $r$-dimensional bar framework $(G,p)$ on $n \geq r+2$ nodes
in general position in $\Rs^r$ admits a positive semidefinite stress matrix with rank
$n-r-1$, then $(G,p)$ is universally rigid. In this paper, we generalize this result in two directions.
First, we extend this result to tensegrity frameworks. Second, we replace the general position assumption
by the weaker assumption that in configuration $p$, each point and its neighbors in $G$ affinely span $\Rs^r$.
\end{abstract}

\section{Introduction}
A {\em tensegrity graph} $G=(V,E)$ is a simple connected graph where $V=\{1,\ldots,n\}$, and
where each edge in $E$ is labelled as either a bar, a cable or a strut.
A {\em tensegrity framework} in $\Rs^r$, denoted by $(G,p)$,
is a tensegrity graph $G=(V, E)$ where each node $i$ is mapped to  a point $p^i$ in
$\Rs^r$. The points $p^1,\ldots,p^n$ will be referred to collectively as the configuration
$p$ of $(G,p)$. $(G,p)$ is {\em $r$-dimensional} if $r$ is the dimension of the affine span of
configuration $p$.
Let $B$, $C$ and $S$ denote the sets of bars, cables and struts of
$(G,p)$ respectively. Then
\[ E= B \cup C \cup S. \]
A {\em bar framework} $(G,p)$ is a tensegrity framework where $E=B$, i.e., $C=S=\emptyset$.

Let $(G,p)$ and $(G,q)$ be two $r$-dimensional and $s$-dimensional tensegrity
frameworks in $\Rs^r$ and $\Rs^s$ respectively. Then
$(G,q)$ is said to be {\em congruent} to $(G,p)$ if:
\[
 ||q^i-q^j||^2 = ||p^i - p^j ||^2 \quad \text{for all $i,j=1,\ldots,n$}.
\]
On the other hand, we say that
 $(G,q)$ is {\em dominated} by $(G,p)$, (or $(G,p)$ {\em dominates} $(G,q)$), if:
\beq \label{defd}
\begin{array}{lcl}
||q^i-q^j||^2 &=& ||p^i - p^j ||^2 \quad \text{for each $\{i,j\} \in B$}, \\
||q^i-q^j||^2 & \leq & ||p^i - p^j ||^2 \quad \text{for each $\{i,j\} \in C$},
\\
||q^i-q^j||^2 &\geq & ||p^i - p^j ||^2 \quad \text{for each $\{i,j\} \in S$}.
\\
\end{array}
\eeq
Moreover, $(G,q)$ is said to be {\em affinely-dominated} by
$(G,p)$ if $(G,q)$ is dominated by $(G,p)$ and $q^i=Ap^i + b$ for all $i=1,\ldots,n$,
where $A$ is an $r \times r$ matrix and $b \in \Rs^r$.

An $r$-dimensional tensegrity framework $(G,p)$ is said to be {\em dimensionally rigid} if no
$s$-dimensional tensegrity framework $(G,q)$, for any $s \geq r+1$, is dominated by $(G,p)$.
Furthermore, $(G,p)$ is said to be {\em universally rigid} if
every $s$-dimensional tensegrity framework $(G,q)$, for any $s$, that is dominated by
$(G,p)$, is in fact congruent to $(G,p)$.

An {\em equilibrium stress} (or simply a stress) of $(G,p)$ is a real-valued function
$\omega$ on $E$, the edge set of $G$, such that
\beq \label{defw}
\sum_{j:\{i,j\} \in E} \omega_{ij} (p^i - p^j) = \bz \mbox{ for all } i=1,\ldots,n.
\eeq

A stress $\omega=(\omega_{ij})$ is said to be {\em proper} if
$\omega_{ij} \geq 0$ for all $\{i,j\} \in C$ and
$\omega_{ij} \leq 0$ for all $\{i,j\} \in S$.
Let $\overline{E}$ denote the set of missing edges in $G$, i.e.,
\[
\overline{E}= \{ \{i,j\}: i \neq j , \{i,j\} \not \in E \},
\]
and let $\omega =(\omega_{ij})$ be a stress of $(G,p)$. Then the
$n \times n$ symmetric matrix $\Omega$ where
\beq \label{defO}
\Omega_{ij} = \left\{ \begin{array}{ll} -\omega_{ij} & \mbox{if } \{i,j\} \in E, \\
                        0   & \mbox{if }  \{i,j\}  \in \overline{E}, \\
                   {\displaystyle \sum_{k:\{i,k\} \in E} \omega_{ik}} & \mbox{if } i=j,
                   \end{array} \right.
\eeq
is called the {\em stress matrix} associated with $\omega$, or a stress matrix
of $(G,p)$. A stress matrix $\Omega$ is {\em proper} if it is associated with a proper
stress $\omega$.
An $n \times n$ matrix $A$ is said to be {\em positive semidefinite} if $x^TA x \geq 0$
for all $x \in \Rs^n$. Furthermore, $A$ is said to be {\em positive definite} if
$x^TAx > 0 $ for all non-zero $x \in \Rs^n$.

The following result provides a sufficient condition for the
universal rigidity of a given tensegrity framework.

\begin{thm}[Connelly \cite{con82}]
\label{thmcon1}
Let $(G,p)$ be an $r$-dimensional tensegrity framework on $n$ vertices in $\Rs^r$, for some $r
\leq n-2$. If the
following two conditions hold:
\begin{enumerate}
\item There exists a proper positive semidefinite stress matrix $\Omega$ of $(G,p)$ of rank
$n-r-1$.
\item There does not exist a tensegrity framework $(G,q)$ in $\Rs^r$
that is affinely-dominated by, but not congruent to, $(G,p)$.
\end{enumerate}
Then $(G,p)$ is universally rigid.
\end{thm}

Connelly \cite{con05} (see also Laurent and Varvitsiotis \cite{lv13})
proved that under the assumption that configuration $p$ is generic, 
Condition 1 of Theorem \ref{thmcon1} implies  Condition 2.
A configuration $p$ is \emph{generic} if the
coordinates of the points $p^1,\ldots,p^n$ are algebraically independent over
the rationals. Thus, for generic tensegrity frameworks, Condition 1 of Theorem \ref{thmcon1} is
sufficient for universal rigidity. When restricted to bar frameworks, i.e., tensegrity frameworks
with no cables or struts, this result
was strengthened, recently, by Alfakih and Ye \cite{ay13}
who proved that Condition 2 of Theorem \ref{thmcon1} is implied by Condition 1 under the weaker assumption
that configuration $p$ is in general position. A configuration $p$ is in general position in $\Rs^r$
for some $r \leq n-1$,
if every subset of $\{p^1,\ldots, p^n\}$ of cardinality $r+1$ is affinely
independent, i.e., every $r+1$ of the points $p^1,\ldots,p^n$ affinely span $\Rs^r$.

\begin{thm}[Alfakih and Ye \cite{ay13}]
\label{thmay}
Let $(G,p)$ be an $r$-dimensional bar framework on $n$ vertices in $\Rs^r$, for some $r
\leq n-2$. If the
following two conditions hold:
\begin{enumerate}
\item There exists a positive semidefinite stress matrix $\Omega$ of $(G,p)$ of rank
$n-r-1$.
\item The configuration $p$ is in general position.
\end{enumerate}
Then $(G,p)$ is universally rigid.
\end{thm}

In this paper, we present characterizations
of dominated and affinely-dominated tensegrity frameworks. As a result,
we strengthen the results in \cite{ay13}
in two directions. First, we extend these results to
tensegrity frameworks. Second, we replace  the general position assumption with the weaker assumption that
in configuration $p$, each point and its neighbors in $G$  affinely span $\Rs^r$.

The paper is organized as follows. In Section \ref{secmain}, we present the main theorems of the paper.
In Section \ref{secpre}, we present the necessary mathematical preliminaries.
The characterizations of dominated and affinely-dominated tensegrity frameworks are given
in Sections \ref{secd} and \ref{secad} respectively. Finally, in Section \ref{secproof}
we present the proofs of the main theorems.

\subsection{Notation}
For easy reference, the notation used throughout the paper are collected below.
$\Sn$ denotes the space of $n \times n$ symmetric matrices equipped with the usual
inner product $\tl A, B \tr$ = trace$\,(AB)$. We denote by $e$ the vector of all
1's in $\Rs^n$, and by $e^i$ the $i$th standard unit vector in $\Rs^n$.
For $i < j$, $F^{ij}= (e^i-e^j)(e^i-e^j)^T$ and $E^{ij}=e^i (e^j)^T + e^j (e^i)^T$.
Also, $L^i= e^i e^T + e (e^i)^T$.
$\bz$ denotes the zero matrix or vector of the appropriate dimension. We use
$|A|$ to denote the cardinality of a finite set $A$. For a node $i$,
$N(i)$ denotes the set of neighbors of $i$, i.e., $N(i)=\{ j : \{i,j\} \in E \}$.
For a tensegrity framework $(G,p)$ with a proper stress matrix $\Omega$,
$C^* =\{ \{i,j\} \in C : \omega_{ij} \neq 0\}$ and
$S^*=\{ \{i,j\} \in S : \omega_{ij} \neq 0\}$. Moreover,
$C^0=\{ \{i,j\} \in C : \omega_{ij} = 0\}$ and
$S^0=\{ \{i,j\} \in S : \omega_{ij} = 0\}$.
The set of missing edges of $G$ is denoted by $\overline{E}$.
$\Es(y) = \sum_{ \{i,j\} \in \overline{E} \cup C \cup S} y_{ij} E^{ij}$ and
$\Es^0(y) = \sum_{ \{i,j\} \in \overline{E} \cup C^0 \cup S^0} y_{ij} E^{ij}$.
Finally, the set theoretic difference is denoted by ``$\backslash$".

\section{Main Results}
\label{secmain}

The following theorem is our main result.

\begin{thm}
\label{thmmain}
Let $(G,p)$ be an $r$-dimensional tensegrity framework on $n$ vertices in $\Rs^r$, for some $r\leq n-2$. 
Then $(G,p)$ is universally rigid
if the following two conditions hold.
\begin{enumerate}
 \item $(G,p)$ admits a proper positive semidefinite stress matrix $\Omega$ with rank $n-r-1$.
 \item For each vertex $i$, the set
 $\{p^i\} \cup \{p^j: \{i,j\}  \in B \cup C^* \cup S^* \}$
affinely spans $\Rs^r$,
\end{enumerate}
where $C^*=\{\{i,j\} \in C : \omega_{ij} \neq 0\}$ and
$S^*=\{\{i,j\} \in S : \omega_{ij} \neq 0\}$.
\end{thm}

The next weaker result is a corollary of Theorem \ref{thmmain}.
\begin{cor}
\label{cormain}
Let $(G,p)$ be an $r$-dimensional tensegrity framework on $n$ vertices in $\Rs^r$, for some $r\leq n-2$. 
Then $(G,p)$ is universally rigid
if the following two conditions hold.
\begin{enumerate}
 \item $(G,p)$ admits a proper positive semidefinite stress matrix $\Omega$ with rank $n-r-1$.
 \item For each vertex $i$, the set $\{p^i\} \cup \{p^j: \{i,j\}  \in B \cup C^* \cup S^* \}$
is in general position in $\Rs^r$.
\end{enumerate}
where $C^*=\{\{i,j\} \in C : \omega_{ij} \neq 0\}$ and
$S^*=\{\{i,j\} \in S : \omega_{ij} \neq 0\}$.
\end{cor}

The proofs of Theorem \ref{thmmain} and Corollary \ref{cormain} are given in Section \ref{secproof}.
In case of bar frameworks, i.e., tensegrity frameworks with no cables or struts,
Theorem \ref{thmmain} and Corollary \ref{cormain} reduce to the following.

\begin{thm}
\label{thmmainb}
Let $(G,p)$ be an $r$-dimensional bar framework on $n$ vertices in $\Rs^r$, for some $r\leq n-2$. 
Then $(G,p)$ is universally rigid
if the following two conditions hold.
\begin{enumerate}
 \item $(G,p)$ has a positive semidefinite stress matrix $\Omega$ with rank $n-r-1$.
 \item For each vertex $i$, the set $\{p^i\} \cup \{p^j: j \in N(i)\}$
affinely spans $\Rs^r$,
\end{enumerate}
where $N(i)$ is the set of adjacent nodes of $i$.
\end{thm}

\begin{cor}
\label{cormainb}
Let $(G,p)$ be an $r$-dimensional bar framework on $n$ vertices in $\Rs^r$, for some $r\leq n-2$. 
Then $(G,p)$ is universally rigid
if the following two conditions hold.
\begin{enumerate}
 \item $(G,p)$ has a positive semidefinite stress matrix $\Omega$ with rank $n-r-1$.
 \item For each vertex $i$, the set $\{p^i\} \cup \{p^j: j \in  N(i)\}$ is in general position in $\Rs^r$.
\end{enumerate}
where $N(i)$ is the set of adjacent nodes of $i$.
\end{cor}

Gortler and Thurston  \cite{gt09} proved that if an $r$-dimensional generic
bar framework $(G,p)$ in $\Rs^r$ is universally rigid,
then $(G,p)$ admits a positive semidefinite stress matrix $\Omega$ of rank $n-r-1$.
An obvious question is to ask whether this result continues to hold if the assumption
of generic configuration is replaced by the assumption of configuration
in general position, i.e., whether the converse of Theorem \ref{thmay} holds.
To this end, Alfakih {\em et al} \cite{aty13} and Alfakih \cite{alf12} proved that indeed
the converse of Theorem \ref{thmay} holds if $G$ is
an $(r+1)$-lateration graph or a chordal graph. Unfortunately,
the converse of Theorem \ref{thmay} does not hold for general graphs
as the following example,
due to Connelly and Whiteley, shows.

\begin{figure}[t]
\setlength{\unitlength}{0.8cm}
\begin{picture}(4,7.5)(-8,-4)
\put(-5,-2){\circle*{0.2}}
\put(-1,0){\circle*{0.2}}
\put( -3,2){\circle*{0.2}}
\put(-4, -1){\circle*{0.2}}
\put(-2, 0){\circle*{0.2}}
\put(-3, 1){\circle*{0.2}}

\put(-5,-2){\line(2,1){4}}
\put(-5,-2){\line(1,1){1}}
\put(-5,-2){\line(1,2){2}}
\put(-4,-1){\line(2,1){2}}
\put(-4,-1){\line(1,2){1}}
\put(-3,1){\line(0,1){1}}
\put(-3,1){\line(1,-1){1}}
\put(-3,2){\line(1,-1){2}}
\put(-2,0){\line(1,0){1}}

\put(-5.4,-2.4){$1$}
\put(-0.6,0){$2$}
\put(-3.2,2.4){$3$}
\put(-4.3,-0.9){$4$}
\put(-2.1,0.2){$5$}
\put(-3.4,1){$6$}

\put(1,-2){\circle*{0.2}}
\put(5,0){\circle*{0.2}}
\put(3,2){\circle*{0.2}}
\put(2, -1){\circle*{0.2}}
\put(4, 0){\circle*{0.2}}
\put(3, 1){\circle*{0.2}}
\put(3.45, 1){\circle*{0.2}}

\put(1,-2){\line(2,1){4}}
\put(1,-2){\line(1,1){1}}
\put(1,-2){\line(1,2){2}}
\put(2,-1){\line(2,1){2}}
\put(2,-1){\line(1,2){1}}
\put(3,1){\line(0,1){1}}
\put(3,1){\line(1,-1){1}}
\put(3,2){\line(1,-1){2}}
\put(4,0){\line(1,0){1}}

\put(3,2){\line(1,-2){0.5}}
\put(3,1){\line(1,0){0.5}}
\put(3.5,1){\line(3,-2){1.5}}

\put(0.6,-2.4){$1$}
\put(5.4,0){$2$}
\put(2.8,2.4){$3$}
\put(1.7,-0.9){$4$}
\put(3.9,0.2){$5$}
\put(2.6,1){$6$}
\put(3.7,1.1){$7$}

\put(-2.9,-3){$(a$)}
\put(3.1,-3){$(b$)}
\end{picture}
\caption{Bar frameworks of Example \ref{cw}. Framework $(a)$ is universally rigid with
    a positive semidefinite stress matrix of rank 3. Framework $(b)$ is also universally
    rigid, however, it does not have a positive semidefinite stress matrix of rank 4. }
\label{fcw}
\end{figure}
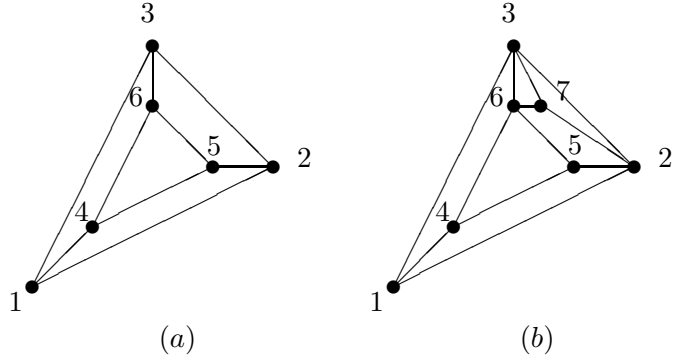
\begin{exa} \label{cw}
Consider the bar framework $(G,p)$ in Figure \ref{fcw}a where
\[p^1=\left[ \begin{array}{r} -2 \\ -2 \end{array} \right], 
p^2=\left[ \begin{array}{r} 2 \\ 0 \end{array} \right], 
p^3=\left[ \begin{array}{r} 0 \\ 2 \end{array} \right], 
p^4=\left[ \begin{array}{r} -1 \\ -1 \end{array} \right], 
p^5=\left[ \begin{array}{r} 1 \\ 0 \end{array} \right], 
\mbox{ and }
p^6=\left[ \begin{array}{r} 0 \\ 1 \end{array} \right]. 
\]

It is easy to check that $(G,p)$ has a unique, up to a scalar multiple, stress matrix
\[
\Omega =\left[ \begin{array}{rrrrrr}
4&1&1&-6&0&0\\
1&4&1&0&-6&0\\
1&1&4&0&0&-6\\
-6&0&0&10&-2&-2\\
0&-6&0&-2&10&-2\\
0&0&-6&-2&-2&10
\end{array} \right].
\]

Note that $\Omega =
\left[ \begin{array}{cc} 3I + ee^T & - 6I \\ -6I & 12 I - 2 ee^T\end{array} \right]$,
where $I$ is the identity matrix of order 3 and $e$ is the vector of all 1's in $\Rs^3$.
Therefore, it follows from Schur's complement \cite{wsv00} that $\Omega$ is positive semidefinite
since $3I+ee^T$ is positive definite and
$(12I-2ee^T) - 6I (3I+ee^T)^{-1} 6 I = \bz$. Moreover, rank $\Omega$ = 3 since $3I+ee^T$
is non-singular. Hence, it follows from Theorem \ref{thmay} that $(G,p)$ is universally rigid.

Now extend $(G,p)$ to framework $(G',p')$ by adding a new
node $7$ adjacent to nodes $2$, $3$ and $6$ such that $(G',p')$ is in general
position as in Figure \ref{fcw}b. Then since $(G,p)$ is universally
rigid, $(G',p')$ must also be universally rigid. Let $\omega'$ be a
non-zero equilibrium stress for $(G',p')$. Then one can see that $\omega_{13}'$ must be
non-zero. So we can assume that  $\omega_{13}'=-1=\omega_{13}$. However, it follows then
that $\omega_{12}'=\omega_{12}$ and $\omega_{14}'=\omega_{14}$.
Hence, $\omega_{45}'=\omega_{45}$ and
$\omega_{46}'=\omega_{46}$. Then again, $\omega_{52}'=\omega_{52}$ and
$\omega_{56}'=\omega_{56}$.
Now, at node 6, since $\omega_{56}'=\omega_{56}$ and
$\omega_{46}'=\omega_{46}$ it follows that
$\omega_{36}'=\omega_{36}$ and $\omega_{67}'= 0$. By repeating this argument
on nodes 3 and 2, we deduce that $\omega_{37}'=\omega_{27}'=0$.
Therefore, the only non-zero equilibrium stress of $(G',p')$
has stress $0$ on every edge incident to node $7$. Therefore, the
corresponding stress matrix $\Omega'$ is obtained from $\Omega$ by adding an all-zero
column and an all-zero row. Hence, obviously rank $\Omega'$ = rank $\Omega$ = $3 < 7-2-1$.
\end{exa}

\section{Preliminaries}
\label{secpre}

In this section we present  mathematical preliminaries that are needed in the sequel.
Let $e^i$ denote the $i$th unit vector in $\Rs^n$ and let
$e$ denote the vector of all 1's in $\Rs^n$. Define the following $n \times n$ symmetric matrices
for $i < j$.
\beq \label{defFEL}
\begin{array}{lcl}
F^{ij} & = & (e^i-e^j)(e^i-e^j)^T,  \\
E^{ij}  & = & e^i (e^j)^T + e^j (e^i)^T, \\
 L^{i}  & = & e^i e^T + e (e^i)^T.
\end{array}
\eeq
Recall that Kronecker delta $\delta_{ij}$ is defined by
\[
\delta_{ij} = \left\{ \begin{array}{ll} 1 & \mbox{if } i=j , \\
                                        0 & \mbox{otherwise}.
                  \end{array} \right.
\]
Then the following two technical lemmas easily follow.
\begin{lem} \label{lem:inner-product}
We have the following equalities.
\begin{equation}
 \trace (F^{kl} E^{ij})  = - 2 \delta_{ki} \delta_{lj},
\end{equation}
for $1\leq i<j\leq n$, $1\leq k< l\leq n$, and
\begin{equation}
\trace (F^{kl} L^{i}) =  0,
\end{equation}
for $1\leq k< l\leq n$ and $1\leq i\leq n$.
\end{lem}
The first equality of Lemma \ref{lem:inner-product} implies that $F^{kl}$ is
orthogonal to $E^{ij}$ if $\{k,l\} \neq \{i,j\}$; while the second equality implies that
$F^{kl}$ is orthogonal to $L^i$ for every $1\leq k<l\leq n$ and $i=1, \ldots, n$.

\begin{lem}{\ } \label{lem:independent}
\begin{enumerate}
\item The set $\{F^{ij} : 1\leq i < j\leq n \}$ is linearly independent.
\item The set $\{E^{ij} : 1\leq i < j\leq n \}\cup
\{L^{i} : 1\leq i\leq n \}$ is linearly independent.
\end{enumerate}
\end{lem}
The following is an immediate consequence of Lemmas \ref{lem:inner-product} and
\ref{lem:independent}.

\begin{cor} \label{corL}
Let $(K, \overline{K})$ be a partition of  the set $ \{\{i,j\} : 1\leq i < j \leq n \}$
and let ${\cal L}$ = span $\{ F^{ij} : \{i,j\} \in K \}$. Then
$\{E^{ij}: \{i,j\} \in \overline{K} \}\cup \{L^i: i=1,\ldots, n\}$
is a basis for ${\cal L}^{\perp}$, the orthogonal complement of
${\cal L}$ in the space of $n \times n$ symmetric matrices $\Sn$.
\end{cor}

\subsection{Stress Matrices and Gale Matrices}

The \emph{configuration matrix} $P$ of  an $r$-dimensional tensegrity framework
$(G,p)$ in $\Rs^r$ is the $n\times r$ matrix whose $i$th row is $(p^i)^T$, i.e.,
\[
P=\left[ \begin{array}{c} (p^1)^T \\ \vdots \\ (p^n)^T \end{array} \right].
\]
Thus the Gram matrix of $(G,p)$ is  $PP^T$.
Moreover,
\begin{equation}
 \|p^i-p^j\|^2 = \trace (F^{ij} PP^T) \quad \text{for $1\leq i<j\leq n$.}
 \label{eq:norm}
\end{equation}
 Without loss of generality we assume the following.
\begin{ass} \label{ass1}
In any configuration $p$, the centroid of the points $p^1,\ldots,p^n$ 
coincides with the origin, i.e., $P^Te= \bz$.
\end{ass}

Since $(G,p)$ is  $r$-dimensional, it follows that $[P \;\; e]$ has full column rank, i.e.,
rank $[P \;\; e] = r+1$. Moreover,
$(G,p)$ is in general position in $\Rs^r$ if and only if every $(r+1) \times (r+1)$
submatrix of $[P \;\; e]$ is non-singular. Furthermore, by definition, we have
that an $n \times n$ matrix $\Omega$ is a stress matrix of $(G,p)$ if
and only if $\Omega_{ij}= 0$ for each $\{i,j\} \in \overline{E}$ and
\[
 \begin{array}{c} P^T \Omega \\ e^T \Omega \end{array}  \begin{array}{l}= \bz, \\ =0. \end{array}
\]
Moreover, it is easy to see from the definition of matrices $F^{ij}$ that
\[
\Omega = \sum_{\{i,j\} \in E} \omega_{ij} F^{ij}.
\]
The dimension of the null space of $\left[ \begin{array}{c} P^T \\ e^T  \end{array} \right]$,
i.e., $n-r-1$, plays an important role. Thus throughout this paper let
\[
\rr = n- r -1.
\]
\begin{defn}
A {\em Gale matrix} of $(G,p)$ is any $n \times \rr$ matrix $Z$ whose columns form
a basis for the null space of $\left[ \begin{array}{c} P^T \\ e^T  \end{array} \right]$.
\end{defn}
\begin{defn}[\cite{gal56,gru67}]
Let $Z$ be a Gale matrix of $(G,p)$ and let ${(z^i)}^T$ be the $i$th row
of $Z$, then $z^i \in \Rs^{\rr}$ is called a {\em Gale transform} of $p^i$
\end{defn}
Hence, we have the following simple but important result.

\begin{lem} \label{lemOZ}
Let $\Omega$ be a stress matrix for an $r$-dimensional framework $(G,p)$, then
\begin{eqnarray*}
&&\Omega = Z \Psi Z^T \mbox{ for some $\rr \times \rr$ symmetric matrix } \Psi,\\
&& \;\; \mbox{ and rank $\Omega$ = rank } \Psi \leq \rr.
\end{eqnarray*}
\end{lem}
Note that Gale matrix is not unique.
In fact, if $Z$ is a Gale matrix of a configuration $p$ and $Q$ is
any non-singular $\rb \times \rb$ matrix then $Z'=ZQ$ is also a Gale matrix of
$p$. Moreover, if $Z,Z'$ are Gale matrices of $p$ then there
exists a non-singular $\rb\times\rb$ matrix $Q$ such that $Z'=ZQ$.

Since Gale matrix $Z$ encodes the affine dependencies among the points
$p^1,\ldots,p^n$, $Z$ has nice properties when some or all of these points
are in general position. The following
lemma is crucial in the proof of our main theorems.

\begin{lem}
\label{lem:spans_vs_nonsingular}
Let $(G,p)$ be an $r$-dimensional tensegrity framework in $\Rs^r$ and let $z^i$ be a Gale transform of
$p^i$ for $i=1,\ldots,n$. Let $J \subset \{1,\ldots,n\}$ such that $|J|=r+1$, and let the set of vectors
$\{p^i: i \in J \}$ be affinely independent. Then the set $\{z^i: i \in {\bar{J}} \}$ is linearly
independent, where $\bar{J}=\{1,\ldots,n\} \setminus J$.
\end{lem}
\bpr
Since $\{p^i: i \in J \}$ is affinely independent, it follows that $[P \;\; e ]_{J}$,
the $(r+1) \times (r+1)$ submatrix of $[P \;\; e]$ whose rows are indexed by $J$, is non-singular.
Let $\lambda_i$, for $i \in \bar{J}$, be scalars such that
$\sum_{i \in \bar{J} }\lambda_i z^i= \bz$. We will show that $\lambda_i=0$ for all $i \in
\bar{J}$. To this end,
set $\lambda_i=0$ for all $i\in J$  and let $\lambda =[\lambda_1 \; \ldots \; \lambda_n]^T$.
Then, by construction, $Z^T \lambda=\bz$. Since the columns of
$[P \;\; e]$ form a basis of the null space of $Z^T$, it follows that
$\lambda= P x + x_0 e$ for some $x\in \Rs^{r}$ and some $x_0 \in \Rs$.
However, by definition of $\lambda$ we have
$(Px+x_0 e)_i = \lambda_i = 0$ for all $i \in J$. Thus,
$x=\bz$ and $x_0=0$ and hence $\lambda = \bz$ and the result follows.
\epr

\begin{cor}
\label{cor:span_vs_ind}
 Let $(G,p)$ be an $r$-dimensional tensegrity framework in $\Rs^r$ and let $z^i$ be a Gale transform of
$p^i$ for $i=1,\ldots,n$. Let $J \subseteq \{1,\ldots,n\}$ and assume that the set of vectors
$\{p^i: i \in J \}$ affinely spans $\Rs^r$. Then the set $\{z^i: i \in \bar{J} \}$
is linearly independent, where $\bar{J} = \{1,\ldots,n\} \setminus J$.
\end{cor}
\bpr
$|J| \geq r+1$ since $\{p^i: i \in J \}$ affinely spans $\Rs^r$.
Let $J'\subseteq J$ such that $|J'|=r+1$ and the set $\{p^i: i\in J'\}$ is affinely
independent. Then by Lemma \ref{lem:spans_vs_nonsingular}, 
the set $\{z^i: i \in \bar{J'}\}$ is linearly independent, where
$\bar{J'} = \{1,\ldots,n\} \backslash J'$. Therefore,
$\{z^i: i \in \bar{J}\}$ is linearly independent since $\bar{J} \subset \bar{J'}$.
\epr

\section{Dominated Tensegrity Frameworks}
\label{secd}
In this section we present
a characterization of tensegrity frameworks dominated by a given $r$-dimensional tensegrity
framework $(G,p)$. Such a characterization, based on Gram matrices,
leads to sufficient conditions for the dimensional and universal rigidities
of $(G,p)$. 

Let $(G,p')$ be an $r'$-dimensional tensegrity framework and let
$P$ and $P'$ be the configuration matrices of $(G,p)$ and $(G,p')$.
For ease of notation, define the following $n \times n$ symmetric matrix

\beq \label{defEs}
\Es(y) = \sum_{\{i,j\} \in \overline{E} \cup C \cup S} y_{ij} E^{ij},
\eeq
 where
$y=(y_{ij}) \in \Rs^{|\overline{E}|+ | C|+ |S|}$.

\begin{lem} \label{lembasic}
Let $P$ and $P'$ be the configuration matrices of tensegrity frameworks $(G,p)$ and $(G,p')$. Then
$(G,p')$ is dominated by $(G,p)$ if and only if
\beq \label{eqG}
P' {P'}^T -PP^T = \Es(y) + x e^T + e x^T,
\eeq
for some  $y=(y_{ij}) \in \Rs^{|\overline{E}|+|C|+|S|}$ and 
$x=(x_i) \in \Rs^n$ where $y_{ij} \geq 0 $ for all $\{i,j\} \in C$ and
$y_{ij} \leq 0 $ for all $\{i,j\} \in S$ and
\beq \label{defx}
x = -\frac{1}{n} \Es(y)e +\frac{1}{2n^2} (e^T \Es(y) e) \; e.
\eeq
\end{lem}

\bpr
Let ${\cal L}$ = span  $\{F^{ij}: \{i,j\} \in B \}$, then  it follows from
Corollary  \ref{corL} that
$\{E^{ij}: \{i,j \} \in \overline{E} \cup C \cup S\} \cup \{L^i: i=1,\ldots,n\}$
is a basis for ${\cal L}^{\perp}$, the orthogonal complement of ${\cal L}$ in $\Sn$.
Now since
\[ ||{p'}^i-{p'}^j||^2 - ||p^i-p^j||^2 = \trace (F^{ij} (P'{P'}^T-PP^T)),
\]
it follows that $(G,p')$ is dominated by $(G,p)$ if and only if
\begin{eqnarray}
 \trace(F^{ij}(P'{P'}^T-PP^T))= 0 \quad \text{for all $\{i,j\} \in B$,} \label{eb} \\
\trace(F^{ij}(P'{P'}^T -PP^T))\leq 0 \quad \text{for all $\{i,j\} \in C$,} \label{ec}\\
\trace(F^{ij}(P'{P'}^T-PP^T)) \geq 0 \quad \text{for all $\{i,j\} \in S$.} \label{es}
\end{eqnarray}
But (\ref{eb}) holds if and only if $(P'{P'}^T-PP^T) \in {\cal L}^{\perp}$, i.e.,
\[
P'{P'}^T-PP^T= \sum_{\{i,j\} \in \overline{E} \cup C \cup S} y_{ij} E^{ij} + \sum_{i=1}^n x_i L^i =
\Es(y) + x e^T + e x^T
\]
 for some $y$ and $x$. Moreover, $(P'{P'}^T-PP^T)e=\bz$
by Assumption \ref{ass1}. Thus
$\Es(y)e + n x + (e^Tx) \; e = 0$. Hence, $(nI+ee^T)x = - \Es(y)e$. Therefore,
\[ x=  (I/n - ee^T/(2n^2)) (-\Es(y) e) =
 -\frac{1}{n} \Es(y)e +\frac{1}{2n^2} (e^T \Es(y) e) \; e.
 \]
Moreover,
(\ref{ec}) holds if and only if $y_{ij} \geq 0$ for all $\{i,j\} \in C$ since
for $\{k,l\} \in C$ it follows from Lemma \ref{lem:inner-product} that
\[
\trace (F^{kl}(P'{P'}^T-PP^T)) = \sum_{\{i,j\} \in C} \trace (F^{kl}E^{ij}) y_{ij} =
-2 y_{kl} \leq 0.
\]
Similarly, (\ref{es}) holds if and only if $y_{ij} \leq 0$ for all $\{i,j\} \in S$.
\epr

The following theorem extends a similar one for bar frameworks \cite{alf11} to tensegrity frameworks.
\begin{thm} \label{th1}
Let $(G,p)$ be a given tensegrity framework and let
$\Omega$ be a proper positive semidefinite stress matrix of $(G,p)$. Then
$\Omega$ is a proper stress matrix for all tensegrity frameworks $(G,p')$ dominated by $(G,p)$.
\end{thm}

\bpr
Let $(G,p')$ be a tensegrity framework dominated by $(G,p)$ and let
$P$ and $P'$ be the Gram matrices of $(G,p)$ and $(G,p')$ respectively.
Then
\begin{eqnarray*}
\omega_{ij} \; \trace(F^{ij} (P'{P'}^T- PP^T)) = 0 & & \mbox{ for each } \{i,j\} \in B, \\
\omega_{ij} \; \trace(F^{ij} (P'{P'}^T- PP^T)) \leq 0 & & \mbox{ for each } \{i,j\} \in C, \\
\omega_{ij} \; \trace(F^{ij} (P'{P'}^T- PP^T)) \leq 0 & & \mbox{ for each } \{i,j\} \in S.
\end{eqnarray*}
Therefore,
\[
\sum_{\{i,j\} \in B \cup C \cup S} \omega_{ij} \; \trace (F^{ij} (P'{P'}^T-PP^T)) 
                               = \trace (\Omega (P'{P'}^T-PP^T)) \leq 0.
\]
But $\Omega P=0$. Therefore, trace$\,(\Omega P'{P'}^T) \leq 0$. However, 
both $P'{P'}^T$ and $\Omega$ are positive semidefinite.
Therefore, trace$\,(\Omega P'{P'}^T) = 0$ and hence $\Omega P'= \bz$. Therefore, 
$\Omega$ is a stress matrix of $(G,p')$.
\epr

The following sufficient condition for dimensional rigidity is an immediate corollary of
Theorem \ref{th1}.
\begin{thm} \label{dr}
Let $(G,p)$ be a given $r$-dimensional tensegrity framework in $\Rs^r$ and let
$\Omega$ be a proper positive semidefinite stress matrix of $(G,p)$ of rank $n-r-1$. Then
$(G,p)$ is dimensionally rigid.
\end{thm}

\bpr
Let $(G,p')$ be an $r'$-dimensional tensegrity framework in $\Rs^{r'}$ dominated by $(G,p)$ and
let $P'$ be the configuration matrix of $(G,p')$. Then dim (null space of $\Omega) = r+1$. Moreover,
it follows from Theorem \ref{th1} that $\Omega P'= \bz$. But $e^T P'=\bz$.  Therefore,
rank $P' \leq r$ since $\Omega e= \bz$. Thus $r' \leq r$ and hence, $(G,p)$ is dimensionally rigid.
\epr

Theorem \ref{th1} also provides an immediate proof of Theorem \ref{thmcon1} as follows.
Let $(G,p')$ be an $r'$-dimensional tensegrity framework in $\Rs^r$ dominated by $(G,p)$ 
and let $P'$ be the configuration matrix
of $(G,p')$. By the proof of Theorem \ref{dr} we have that the columns of $P'$ belong to the null space
of $\left[ \begin{array}{c} \Omega \\ e^T \end{array} \right]$. Then there exists an $r \times r$
matrix $A$ such that $P'=PA$ since the columns of $P$ form a basis of this null space.
Hence, configuration $p'$ is obtained from $p$ by an affine transformation. Hence, $(G,p')$ is
congruent to $(G,p)$ and thus $(G,p)$ is universally rigid.
\epr

\section{Affinely-Dominated Tensegrity Frameworks}
\label{secad}

Recall that an affine motion is of the form ${p'}^i=Ap^i+b$ for all $i=1,\ldots,n$.
In this section we characterize tensegrity frameworks that are affinely-dominated by
a given $r$-dimensional tensegrity framework $(G,p)$.
The following lemma characterizes affine-domination in terms of configuration $p$ and
the bars, cables and struts of $(G,p)$.

\begin{lem} \label{lem:affP}
Let $(G,p)$ be an $r$-dimensional tensegrity framework in $\Rs^r$. Then there exists
a tensegrity framework $(G,p')$ affinely-dominated by, but not congruent to, $(G,p)$ if and only if
there exists a non-zero symmetric $r \times r$ matrix $\Phi$ such that
\begin{eqnarray}
&&\trace (F^{ij}(P \Phi P^T)) = 0 \quad \text{for all $\{i,j\} \in B$}, \label{eb2} \\
&&\trace (F^{ij}(P \Phi P^T)) \leq 0 \quad \text{for all $\{i,j\} \in C$}, \label{ec2} \\
&&\trace (F^{ij} (P \Phi P^T)) \geq 0 \quad \text{for all $\{i,j\} \in S$}. \label{es2}
\end{eqnarray}
\end{lem}

\bpr To prove the ``only if" part
assume that $(G,p')$ is affinely-dominated by, but not congruent to, $(G,p)$. Then
$P'=PA$ for some $r \times r$ matrix $A$ since ${P'}^Te=0$. Thus by (\ref{eq:norm})
\[
||{p'}^i -{p'}^j||^2-||p^i - p^j||^2  = \trace (F^{ij}(P (A A^T - I) P^T)).
\]
The result follows by setting $\Phi=AA^T-I$. Obviously,
$\Phi$ is symmetric. Furthermore, $A$ is not orthogonal since $(G,p')$ is not congruent
to $(G,p)$. Thus $\Phi \neq \bz$.

To prove the ``if" part assume that there exists a non-zero symmetric matrix $\Phi$ satisfying
(\ref{eb2})--(\ref{es2}). Then there exists a sufficiently small $\epsilon >0$
such that $I + \epsilon \Phi$ is positive definite. Thus there exists an $r \times r$ nonsingular matrix $A$
such that $AA^T=I+ \epsilon \Phi$. Hence,
\[
\trace (F^{ij}(P \Phi P^T)) = \frac{1}{\epsilon} \, \trace (F^{ij} (P (AA^T-I) P^T)) =
\frac{1}{\epsilon} \, (||{p'}^i -{p'}^j||^2-||p^i - p^j||^2).
\]
Thus the result follows.
\epr

Affine-domination can also be characterized in terms of Gale matrix $Z$ and the missing edges,
cables and struts of $(G,p)$.

\begin{lem} \label{lem:affZ}
Let $(G,p)$ be an $r$-dimensional tensegrity framework in $\Rs^r$ and let $Z$ be a Gale matrix
of $(G,p)$. Then there exists a tensegrity framework
$(G,p')$ affinely-dominated by, but not congruent to, $(G,p)$ if and only if
there exists a non-zero $y=(y_{ij}) \in \Rs^{|\overline{E}|+ | C|+ |S|}$ and
$\xi =(\xi_i) \in \Rs^{\rr}$  where $y_{ij} \geq 0 $ for all $\{i,j\} \in C$ and
$y_{ij} \leq 0 $ for all $\{i,j\} \in S$ such that
\beq \label{eq:affZ}
\Es(y)Z = e \xi^T.
\eeq
\end{lem}

\bpr Assume that there exists a non-zero symmetric matrix $\Phi$ satisfying
(\ref{eb2})--(\ref{es2}).
Let ${\cal L}$ = span  $\{F^{ij}: \{i,j\} \in B \}$. Then  it follows from
Corollary \ref{corL} that
$\{E^{ij}: \{i,j\} \in \overline{E} \cup C \cup S\} \cup \{L^i: i=1,\ldots,n\}$
is a basis for ${\cal L}^{\perp}$, the orthogonal complement of ${\cal L}$ in $\Sn$.
Since $\trace (F^{ij} P \Phi P^T)= 0 $ for all $\{i,j\} \in B$ if and
only if
$P \Phi P^T \in {\cal L}^\perp$, we have
\[
P \Phi P^T= \sum_{(i,j) \in \overline{E} \cup C \cup S} y_{ij} E^{ij} + \sum_{i=1}^n x_i L^i =
\Es(y) + x e^T + e x^T
\]
 for some $y$ and $x$. Moreover, since $P \Phi P^T e= \bz$, it follows that
$\Es(y)e + n x + (e^Tx) \, e = \bz $. Hence,
\[
x = -\frac{1}{n} \Es(y)e +\frac{1}{2n^2} (e^T \Es(y) e) \, e.
 \]
 (see the proof of Lemma \ref{lembasic} ).
 Note that $y \neq \bz$ since $\Phi \neq \bz$.
Furthermore,
\[ \trace \,(F^{ij}P \Phi P^T) =
\trace \,(F^{ij} ( \Es(y) + x e^T + e x^T))  = -2 y_{ij}.
\]
Therefore,  $\trace \,(F^{ij}P \Phi P^T) \leq 0 $ for all $\{i,j\} \in C$ if and only if
$y_{ij} \geq 0 $ for all $\{i,j\} \in C$ and, similarly,
 trace$ \,(F^{ij}P \Phi P^T) \geq 0 $ for all $\{i,j\} \in S$ if and only if
$y_{ij} \leq 0 $ for all $\{i,j\} \in S$.
Therefore, $y$ and $\xi = -Z^Tx$ is a solution of (\ref{eq:affZ}) since
$\Es(y)Z + e \; x^T Z =P \Phi P^TZ = \bz$.

Conversely, assume that there exists a solution of (\ref{eq:affZ}), where $y \neq 0$ and
$y_{ij} \geq 0$ for all $\{i,j\} \in C$ and $y_{ij} \leq 0$ for all $\{i,j\} \in S$.
Then $P^T \Es(y)Z = \bz$ and $Z^T \Es(y) Z = \bz$. Hence,
\beq \label{eq:rho}
\begin{array}{ll}
\Es(y) & =  P \Phi P^T + P \zeta e^T + e \zeta^T P^T + Z \rho e^T + e \rho^T Z^T + \sigma ee^T, \\
       & =  P \Phi P^T + (P \zeta + Z \rho + \sigma e/2 ) e^T + e (P \zeta + Z \rho + \sigma e /2)^T,
\end{array}
\eeq
for some symmetric matrix $\Phi$ and vectors $\zeta$ and $\rho$ and scalar $\sigma$.
But, by multiplying (\ref{eq:rho}) from the right by $e$, we get
$\Es(y)e/n =   P \zeta  +  Z \rho + \sigma e$. Moreover, by
multiplying (\ref{eq:rho}) from the left by $e^T$ and from the right by $e$ we
get
$e^T\Es(y)e = n^2 \sigma$. Therefore,
\[
\Es(y) = P \Phi P^T - x e^T - e x^T,
\]
where $x$ is as given in (\ref{defx}).
Note that $\Phi \neq 0$. Thus $\Phi$ satisfies (\ref{eb2})--(\ref{es2})
and the result follows from Lemma \ref{lem:affP}.
\epr

Lemmas  \ref{lem:affP} and \ref{lem:affZ}
 can be strengthened if a non-zero proper stress matrix $\Omega$ of $(G,p)$ is known,
or if rank $\Omega = n-r-1$. We discuss these two cases in the next two subsections.

\subsection{The Case Where $\Omega$ is Known}

In this subsection we assume that a non-zero proper stress matrix of $(G,p)$ is known.
Then we have the following lemmas.

\begin{lem} \label{lem:affP2}
Let $(G,p)$ be an $r$-dimensional tensegrity framework in $\Rs^r$ and let $\Omega$ be a proper stress
matrix of $(G,p)$. Then there exists a tensegrity framework
$(G,p')$ affinely-dominated by, but not congruent to, $(G,p)$ if and only if
there exists a non-zero symmetric $r \times r$ matrix $\Phi$ such that:
\begin{eqnarray}
&&\trace (F^{ij} P \Phi P^T) = 0 \quad \text{for all $\{i,j\} \in B \cup C^* \cup S^* $,} \label{eb3} \\
&&\trace (F^{ij} P \Phi P^T) \leq 0 \quad \text{for all $\{i,j\} \in C^0$,} \label{ec3} \\
&&\trace (F^{ij} P \Phi P^T) \geq 0  \quad \text{for all $\{i,j\} \in S^0$,}
\label{es3}
\end{eqnarray}
where $C^* =\{ \{i,j\} \in C: \omega_{ij} \neq 0\}$,
$S^* =\{ \{i,j\} \in S: \omega_{ij} \neq 0\}$, $C^0 = C \setminus
C^*$ and $S^0 = S \setminus S^*$.
\end{lem}

\bpr
Assume that $\Phi$ satisfies (\ref{eb2})--(\ref{es2}) in Lemma \ref{lem:affP}.
Then, by definition, $\Omega P\Phi P^T = \bz$. Therefore,
\begin{eqnarray*}
\trace (\Omega P \Phi P^T) & = & \sum_{\{i,j\} \in E}
                           \omega_{ij} \; \trace (F^{ij}P \Phi P^T), \\
 & = & \sum_{\{i,j\} \in C \cup S} \omega_{ij} \; \trace (F^{ij}P \Phi P^T), \\
 & = & 0,
\end{eqnarray*}
since $\trace \,(F^{ij} P \Phi P^T) = 0$ for every $\{i,j\} \in B$.
But $\omega_{ij} \; \trace \, (F^{ij} P\Phi P^T) \leq 0$ for every
$\{i,j\} \in C$ since $\omega_{ij} \geq 0$ and $\trace \,(F^{ij} P \Phi P^T) \leq 0$.
Similarly,
$\omega_{ij} \; \trace \,(F^{ij} P \Phi P^T) \leq 0$ for every $(i,j) \in S$.
Therefore,
\[
  \omega_{ij} \; \trace \,(F^{ij} P \Phi P^T) = 0 \quad \text{for each $\{i,j\} \in C \cup S$.}
\]
Thus, $\trace \, (F^{ij} P\Phi P^T) = 0 \quad \text{for each $\{i,j\} \in C^*\cup S^*$}$,
and the result follows.
\epr

Note that
the necessity of (\ref{eb3}) for the existence of $(G,p')$ affinely-dominated by, but not congruent to,
$(G,p)$ was given in Whiteley \cite{whi87}. Also,
it was implicitly given in Laurent and Varvitsiotis \cite{lv13}. 

\begin{lem} \label{lem:affZ2}
Let $(G,p)$ be an $r$-dimensional tensegrity framework in $\Rs^r$ and let $Z$ be a Gale
matrix of $(G,p)$. Then there exists a tensegrity framework
$(G,p')$ affinely-dominated by, but not congruent to, $(G,p)$ if and only if
there exists a non-zero $y=(y_{ij}) \in \Rs^{|\overline{E}|+ | C^0|+ |S^0|}$ and
$\xi =(\xi_i) \in \Rs^{\rr}$  where $y_{ij} \geq 0 $ for all $\{i,j\} \in C^0$ and
$y_{ij} \leq 0 $ for all $\{i,j\} \in S^0$ such that
\beq \label{eq:affZ2}
\Es^0(y)Z = e \xi^T,
\eeq
where $\Es^0(y) = \sum_{\{i,j\} \in \overline{E} \cup C^0 \cup S^0} y_{ij} E^{ij}$.
\end{lem}
\bpr
The proof is identical to that of Lemma \ref{lem:affZ} where in this case
${\cal L}=\mathrm{span}\{F^{ij}: \{i,j\} \in B \cup C^* \cup S^*\}$.
Thus
$\{E^{ij}: \{i,j\} \in \overline{E} \cup C^0 \cup S^0\} \cup \{L^i: i=1,\ldots,n\}$
is a basis for ${\cal L}^{\perp}$.
\epr

The following example is an illustration of Lemma \ref{lem:affZ2}.
\begin{figure}[t]
\thicklines
\setlength{\unitlength}{0.8cm}
\begin{picture}(3,6)(-8,-4)
\put(-2.5,0){\circle*{0.2}}
\put(0,0){\circle*{0.2}}
\put( 2.5,0){\circle*{0.2}}
\put(0,2.5){\circle*{0.2}}

\multiput(-2.5,0) (0.4,0){7}{\line(1,0){0.2}}
\multiput(0,0) (0.4,0){7}{\line(1,0){0.2}}
\qbezier(-2.5,0)(0,-0.5) (2.5,0)
\qbezier(-2.5,-0.1)(0,-0.6) (2.5,-0.1)
\put(0,0){\line(0,1){2.5}}

\put(-2.5,-0.05) {\line(1,1){2.5}}
\put(-2.5,0.05) {\line(1,1){2.5}}

\put(0,2.45) {\line(1,-1){2.5}}
\put(0,2.55) {\line(1,-1){2.5}}

\put(-2.9,0){$1$}
\put(-0.4,0.2){$2$}
\put(2.7,0){$3$}
\put(-0.3,2.8){$4$}

\end{picture}
\caption{The $2$-dimensional tensegrity framework in $\Rs^2$ of Example \ref{exa1}.
    Bars, cables and struts are drawn, respectively, as solid lines, dashed lines and
    double lines.
    The strut $\{1,3\}$ is shown as an arc to make cables $\{1,2\}$ and $\{2,3\}$ visible. }
\label{f1}
\end{figure}
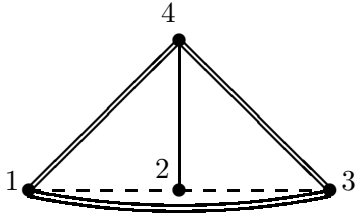

\begin{exa} \label{exa1}

Consider the $2$-dimensional tensegrity framework $(G,p)$ in Figure \ref{f1}.
Then, obviously, there does not exist a tensegrity framework $(G,p')$ affinely-dominated
by $(G,p)$. Next we show how this is implied by Lemma \ref{lem:affZ2}.
In this case the sets of cables and struts with non-zero stresses are given by
$C^* = \{  \{1,2\}, \{2,3\}  \}$ and $S^* = \{  \{1,3\}  \}$. Thus,
$\overline{E}$ = $C^0 = \emptyset$ and $S^0 = \{  \{1,4\}, \{3,4 \}  \}$.
Hence,
\[
\Es^0(y) = y_{14} E^{14}+ y_{34} E^{34} = \left[ \begin{array}{cccc} 0 & 0 & 0 & y_{14} \\ 0 & 0 & 0 & 0 \\
                                   0 & 0 & 0 & y_{34} \\ y_{14} & 0 & y_{34} & 0 \end{array} \right],
\mbox{ and } Z =\left[ \begin{array}{r} 1 \\ -2 \\ 1 \\ 0 \end{array} \right].
\]
Thus, $\Es^0(y)Z=e \xi$ reduces to $y_{14}+y_{34}=0$. Hence, the only solution to $\Es^0(y)Z= e \xi $
where $y_{14} \leq 0$, and $y_{34} \leq 0$ is the trivial solution $y_{14}=y_{34}=0$. Hence,
by Lemma \ref{lem:affZ2}, there does not exist a tensegrity framework $(G,p')$ that is affinely-dominated
by, but not congruent to, $(G,p)$.

Now suppose that strut $\{3,4\}$ is replaced by a cable. Then obviously, in this case,
bar $\{2,4\}$ can rotate to the right.
On the other hand, $\Es^0(y)Z= e \xi $
where $y_{14} \leq 0$, and $y_{34} \geq 0$ has a non-zero solution
where $y_{14} =-1$, and $y_{34}=1$. Hence, by Lemma \ref{lem:affZ2},
there exists a framework $(G,p')$ that is affinely dominated by,
but not congruent to, $(G,p)$.
\end{exa}

\subsection{The Case where rank $\Omega = n-r-1$}

In this subsection we assume that $\Omega$ is a proper stress matrix of $(G,p)$ and
rank $\Omega = \rr$. We begin with the following lemma which establishes the
existence of a special Gale matrix $\Zh$ with desirable properties.

\begin{lem}
\label{lemZh1}
Let $(G,p)$ be an $r$-dimensional tensegrity framework in $\Rs^r$ and let
$\Omega$ be a stress matrix of $(G,p)$ with rank $\rr$.
Then there exists an index set
$J=\{j_1, \ldots, j_{\rr} \} \subset \{1,\ldots,n\}$ and a Gale matrix
$\hat{Z}$ of $(G,p)$ whose columns are indexed by $J$. Furthermore,
$\hat{Z}$ has the following property:
\beq
\mbox{ for } k=1,\ldots,\rr, \;\; \zh_{ij_k}=0 \mbox{ for each } \{i,j_k\} \in \overline{E}.
\eeq
\end{lem}

\bpr
Since rank $\Omega = \rr$, then there exist $\rr$ linearly independent columns of $\Omega$.
Let these columns be indexed by $J=\{j_1,\dots, j_{\rb}\}$ and let $\hat{Z}$
be the $n \times \rr$ submatrix of $\Omega$ whose columns are indexed by $J$.
Then $\Zh$ is a Gale matrix of $(G,p)$
since $P^T \Omega = \bz$ and $e^T \Omega= \bz$, and since $\hat{Z}$ has full column rank.
Furthermore, since $\Omega_{ij}=0$ for $\{i,j\} \in \overline{E}$, it follows
that $\zh_{ij_k}=0$ for each $\{i,j_k\} \in \overline{E}$ and for each $k=1,\ldots,\rr$.
\epr

The following lemma, which is a stronger version of Lemma \ref{lem:affZ2} in case rank $\Omega = \rr$,
is key to our proof of Theorem \ref{thmmain}.

\begin{lem} \label{lem:affZ3}
Let $(G,p)$ be an $r$-dimensional tensegrity framework in $\Rs^r$ and assume that
$\Omega=Z \Psi Z^T$ is a proper stress matrix of $(G,p)$ with rank $n-r-1$.
Then there exists a tensegrity framework
$(G,p')$ affinely-dominated by, but not congruent to, $(G,p)$ if and only if
there exists a non-zero $y=(y_{ij}) \in \Rs^{|\overline{E}|+ | C^0|+ |S^0|}$
where $y_{ij} \geq 0 $ for all $\{i,j\} \in C^0$ and
$y_{ij} \leq 0 $ for all $\{i,j\} \in S^0$ such that
\beq \label{eq:affZ3}
\Es^0(y) Z =\bz,
\eeq
where $\Es^0(y) = \sum_{\{i,j\} \in \overline{E} \cup C^0 \cup S^0} y_{ij} E^{ij}$.
\end{lem}

\bpr
Let $\hat{Z}$ be the Gale matrix of Lemma \ref{lemZh1}. Then it suffices to show that
$\Es^0(y) \hat{Z} = e \xi^T$ is equivalent to  $\Es^0(y) \hat{Z} = \bz$.
It is trivial that if $\Es^0(y) \hat{Z} =\bz $ implies $\Es^0(y) \hat{Z} = e \xi^T$.
Next we prove that
$\Es^0(y) \hat{Z} = e \xi^T$ implies  $\Es^0(y) \hat{Z} = \bz$.
To this end,
for every $k=1,\dots,\rb$ we have
\begin{eqnarray*}
 \xi_{j_k}&=&\sum_{i=1}^n (\Es^0(y))_{j_k i} \; \zh_{ij_k},  \\
&=&      (\Es^0(y))_{j_kj_k} \; \zh_{j_kj_k} \\
   &&    +\sum_{i:\{i,j_k\} \in B \cup C^* \cup S^* } (\Es^0(y))_{j_ki} \; \zh_{ij_k} \\
   &&     + \sum_{i: \{i,j_k\} \in \overline{E} \cup C^0 \cup S^0 } (\Es^0(y))_{j_ki} \; \zh_{ij_k}\\
&=&  0,
\end{eqnarray*}
since $(\Es^0(y))_{j_k j_k} = 0$ and
since if $\{i,j_k\} \in B \cup C^* \cup S^*$, then $(\Es^0(y))_{ij_k}=0$; and if
$\{ i,j_k\} \in \overline{E}$, then
$\zh_{ij_k}=0$. Moreover, if $\{ i,j_k\} \in C^0 \cup S^0$, then
$\zh_{ij_k}=-\omega_{ij_k}=0$.
Therefore $\Es^0(y) \Zh=0$.
\epr

\section{Proofs of Theorem \ref{thmmain} and Corollary \ref{cormain}}
\label{secproof}

\noindent {\bf Proof of Theorem \ref{thmmain}}

Suppose that an $r$-dimensional tensegrity framework $(G,p)$ in $\Rs^r$ admits a positive
semidefinite stress matrix $\Omega$ of rank $\rr$, and suppose that for each node $i$, the set
$\{p^i\} \cup \{p^j: \{i,j\} \in B \cup C^* \cup S^* \}$ affinely spans $\Rs^r$.
Then  it suffices to show that
Condition 2 of Theorem \ref{thmcon1} holds, i.e., it suffices to show that
there is no tensegrity framework $(G,p')$ affinely-dominated by, but not congruent to, $(G,p)$.
However, by Lemma \ref{lem:affZ3}, Condition 2 of Theorem \ref{thmcon1} holds if the only solution of
the equation $\Es^0(y) Z = \bz$ is the trivial solution $y=\bz$.

To this end, the set $\{p^i\} \cup \{p^j: \{i,j\} \in B \cup C^* \cup S^* \}$ 
affinely spans $\Rs^r$ for each
node $i$. Then by Corollary \ref{cor:span_vs_ind}, the set
$\{z^j: \{i,j \} \in \overline{E} \cup C^0 \cup S^0 \} $ is linearly independent
for each $i=1,\ldots,n$.
Now equation $\Es^0 (y) Z=\bz$ can be written as
$\sum_{j=1}^n (\Es^0(y))_{ij} \; z^j= \bz$ for each $i=1,\ldots,n$,
which is equivalent to
$\sum_{j : \{i,j\} \in \overline{E} \cup C^0 \cup S^0} (\Es^0(y))_{ij} \; z^j=\bz$,
since $(\Es^0(y))_{ij}=0$ for all $j: \{i,j\} \in B \cup C^* \cup S^*$. Therefore, the
linear independence of the set $\{z^j: \{i,j\} \in \overline{E} \cup C^0 \cup S^0 \}$ implies that
$y_{ij}=0$ for all $j$ such that $\{i,j\} \in \overline{E} \cup C^0 \cup S^0$ and
for every $i=1,\dots,n$. Therefore $y=\bz$. This completes the proof.
\epr

\noindent{\bf Proof of Corollary \ref{cormain}}

Corollary \ref{cormain} follows from the following lemma.

\begin{lem} \label{lemdeg}
Let $(G,p)$ be an $r$-dimensional tensegrity framework on $n$ vertices in $\Rs^r$, for $r \leq n-2$,
and let $\Omega$ be a stress matrix of $(G,p)$ with rank $n-r-1$. Assume that
for each node $i=1,\ldots,n$,
the set $\{p^i\} \cup \{p^j: \{i,j\} \in B \cup C^* \cup S^* \; \}$ is in general
position in $\Rs^r$. Then the set $\{p^i\} \cup \{p^j: \{i,j\} \in B \cup C^* \cup S^* \; \}$
affinely spans $\Rs^r$.
\end{lem}

\bpr
It suffices to show that the cardinality of the set
$\{i\} \cup \{j: \{i,j\} \in B \cup C^* \cup S^*  \}$
is at least $r+1$ for each node $i$ of $G$. To this end,
let $i$ be a node of $G$ and let
$\{j: \{i,j\} \in B \cup C^* \cup S^* \}$ = $\{j_1, \ldots, j_k\}$. Then,
under the lemma's assumption,
the set $\{p^i\} \cup \{ p^{j_1}, \ldots, p^{j_k} \}$ is affinely dependent only if
$k+1 \geq r+2$. Thus, the set $\{ p^{j_1}-p^i, \ldots, p^{j_k}-p^i \}$
is linearly dependent only if $k \geq r+1$.

Let $\omega=(\omega_{ij})$ be the stress of $(G,p)$ associated with $\Omega$.
Then since 
\[
\sum_{j: \{i,j\} \in B \cup C^* \cup S^*} \omega_{ij} (p^i - p^j) =
\sum_{j \in \{j_1,\ldots,j_k\} } \omega_{ij} (p^i - p^j) = \bz,
\]
it follows that
either $k \geq r+1$, or $\omega_{ij} = 0$ for all $j$ such that $ \{i,j\} \in B \cup C^* \cup S^*$. 
In the first case we are done, so
assume the latter. Then
the entries of the $i$th row and the $i$th column of $\Omega$ are all 0's. Let
$\Omega'$ be the matrix obtained from $\Omega$ by deleting its $i$th row and $i$th
column. Hence, rank $\Omega'$ = rank $\Omega = n-r-1$.

Let $G'=G - i$, i.e., $G'$ is the graph obtained from $G$ by deleing node $i$ and all
the edges incident with it. Also denote the configuration $ \{p^1,\ldots,p^n\} \backslash \{p^i\}$
by $p'$. 
Since $\Omega \neq \bz$, there exists one node, say $v$, such that the cardinality of the set
$\{v\} \cup \{j: \{v,j\} \in B \cup C^* \cup S^* \}$
is at least $r+2$. Therefore, configuration $p'$ also affinely
spans $\Rs^r$. Thus $(G',p')$ is an $r$-dimensional tensegrity framework in $\Rs^r$ and
$\Omega'$ is a stress matrix of $(G',p')$. Therefore,
it follows from Lemma \ref{lemOZ} that $n-r-1$ = rank $\Omega' \leq (n-1) - r -1$,
a contradiction. Thus the result follows.
\epr

\noindent{\bf \large Acknowledgements}

This paper was motivated by a question of Bob Connelly, asking whether it is possible to extend
Theorem \ref{thmay} to tensegrity frameworks, and whether the general position assumption can be
replaced by a weaker one. The authors would like to
thank Bob Connelly and Walter Whiteley for sending to the unpublished manuscript \cite{whi87}.
Also, thanks to Anthony Man-Cho So for communicating to us Example \ref{cw} 
(due to Connelly and Whiteley).


\end{document}